\documentclass{article}
\usepackage{amsmath}
\usepackage{amssymb}
\usepackage{amsthm}
\usepackage[mathlines]{lineno}

\def\R{\ensuremath{{\bf R}}}
\def\Z{\ensuremath{{\bf Z}}}

\def\N{\ensuremath{{\bf N}}}

\def\P{\ensuremath{{\bf P}}}
\def\loc{\ensuremath{\text{loc}}}
\def\dist{\ensuremath{\text{dist}}}
\newcommand{\floor}[1]{{\ensuremath{[#1]}}}

\def\epsilon{\varepsilon}

\def\lcst{\ensuremath{\underline{\Theta}}}
\def\ucst{\ensuremath{\overline{\Theta}}}
\def\lbi{\ensuremath{\underline{s}}}
\def\ubi{\ensuremath{\overline{s}}}
\newtheorem{thm}{Theorem}
\newtheorem{pro}[thm]{Proposition}
\newtheorem{lem}[thm]{Lemma}
\newtheorem{cor}[thm]{Corollary}
\theoremstyle{definition}
\newtheorem{defn}{Definition}
\theoremstyle{remark}
\newtheorem{rem}{Remark}
\newtheorem{exmp}{Example}
\theoremstyle{plain}

\def\BExe{\begin{exmp}}
\def\EExe{\end{exmp}}
\def\BRem{\begin{rem}}
\def\ERem{\end{rem}}
\def\BThe{\begin{thm}}
\def\EThe{\end{thm}}
\def\BDef{\begin{defn}}
\def\EDef{\end{defn}}
\def\BPro{\begin{pro}}
\def\EPro{\end{pro}}
\def\BLem{\begin{lem}}
\def\ELem{\end{lem}}
\def\BCor{\begin{cor}}
\def\ECor{\end{cor}}
\def\BProof{\begin{proof}}
\def\EProof{\end{proof}}
\begin{document}

\author{ \small \sc Damien Kreit\thanks{Universit\'e de Li\`ege, Institut de Math\'ematique, Grande Traverse, 12, B\^atiment B37, B-4000 Li\`ege (Sart-Tilman), Belgium.} and Samuel Nicolay\thanks{Universit\'e de Li\`ege, Institut de Math\'ematique, Grande Traverse, 12, B\^atiment B37, B-4000 Li\`ege (Sart-Tilman), Belgium.} \thanks{Corresponding author. Email: S.Nicolay@ulg.ac.be. Phone: +32(0)43669433. Fax: +32(0)43669547.}}
\title{Generalized pointwise H\"older spaces}
\maketitle 

\begin{abstract}
In this paper, we introduce a generalization of the pointwise H\"older spaces. We give alternative definitions of these spaces, look at their relationship with the wavelets and introduce a notion of generalized H\"older exponent.
\end{abstract}

\section{Introduction}
In \cite{unifhold,unifhold2}, the properties of generalized uniform H\"older spaces have been investigated. The idea underlying the definition is to replace the exponent $\alpha$ of the usual spaces $\Lambda^\alpha(\R^d)$ (see e.g.~\cite{kra:83}) with a sequence $\sigma$ satisfying some conditions. The so-obtained spaces $\Lambda^\sigma(\R^d)$ generalize the spaces $\Lambda^\alpha(\R^d)$; the spaces $\Lambda^\sigma(\R^d)$ are actually the spaces $B^{1/\sigma}_{\infty,\infty}(\R^d)$, but they present specific properties (induced by $L^\infty$-norms) when compared to the more general spaces $B^{1/\sigma}_{p,q}(\R^d)$ studied in \cite{1076.46025,farkas,1108.42007,1116.46024,kuhn2006,1133.46019} for example. Indeed it is shown in \cite{unifhold,unifhold2} that most of the usual properties holding for the spaces $\Lambda^\alpha(\R^d)$ can be transposed to the spaces $\Lambda^\sigma(\R^d)$.

Here, we introduce the pointwise version of these spaces: the spaces $\Lambda^{\sigma,M}(x_0)$, with $x_0 \in\R^d$. Let us recall that a function $f\in L^\infty_\loc(\R^d)$ belongs to the usual pointwise H\"older space $\Lambda^\alpha(x_0)$ ($\alpha>0$) if and only if there exist $C,J>0$ and a polynomial $P$ of degree at most $\alpha$ such that
\begin{equation}\label{eq:us hol}
 \sup_{|h|\le 2^{-j}} |f(x_0+h)-P(h)| \le C 2^{-j\alpha}.
\end{equation}
As in \cite{unifhold,unifhold2}, the idea is again to replace the sequence $(2^{-j\alpha})_j$ appearing in this inequality with a positive sequence $(\sigma_j)_j$ such that $\sigma_{j+1}/\sigma_j$ and  $\sigma_j/\sigma_{j+1}$ are bounded (for any $j$); the number $M$ stands for the maximal degree of the polynomial (this degree can not be induced by a sequence $\sigma$). By doing so, one tries to get a better characterization of the regularity of the studied function $f$; a usual choice is to replace $2^{-j\alpha}$ with $j2^{-j\alpha}$ (see e.g.~\cite{0873.42019,jaffard,marianne}). Generalizations of the pointwise H\"older spaces have already been proposed (see e.g.~\cite{marianne}), but, to our knowledge, the definition we give here is the most general version and leads to the sharpest results.


This work is organized as follows. We first give the definitions leading to generalized pointwise H\"older spaces $\Lambda^{\sigma,M}(x_0)$ and prove that, under some general conditions, the polynomials appearing in the definition are independent from the scale, as it is the case with the usual H\"older spaces. Next we give some alternative definitions of the spaces $\Lambda^{\sigma,M}(x_0)$, mimicking the different possible definitions of $\Lambda^\alpha(x_0)$. One of the nicest properties of the H\"older spaces is their relationship with the wavelet theory given in \cite{jaffard}; we show here that this result still holds in the general case. Finally, we give some conditions under which one gets embedded generalized pointwise H\"older spaces and define a generalized H\"older exponent.

Throughout this paper, $B$ denotes the open unit ball centered at the origin; moreover we set $B_j=2^{-j} B$. The floor function is denoted $\floor{\cdot}$ and $\P[\alpha]$ designates the set of polynomials of degree at most $\floor{\alpha}$. We use the letter $C$ for generic positive constant whose value may be different at each occurrence.

\section{Pointwise generalized H\"older spaces}
To present the generalized pointwise H\"older spaces, we first need to recall some notions concerning the admissible sequences. After having introduced the definitions, we point out a major difference between the usual spaces and the generalized ones: the polynomial arising in the definition depends on the scale. It is then natural to look under which condition this constraint can be dropped.

\subsection{Definition}
The generalization of the H\"older spaces we propose here is based on the notion of admissible sequence \cite{1133.46019}.
\BDef
A sequence $\sigma=(\sigma_j)_{j\in\N}$ of real positive numbers is called admissible if and only if there exists a positive constant $C$ such that
\[
 C^{-1} \sigma_j \le \sigma_{j+1} \le C \sigma_j,
\]
for any $j\in\N$.
\EDef
If $\sigma$ is such a sequence, we set
\[
 \lcst_j=\inf_{k\in\N} \frac{\sigma_{j+k}}{\sigma_k}
 \quad\text{and}\quad
 \ucst_j=\sup_{k\in\N} \frac{\sigma_{j+k}}{\sigma_k}
\]
and define the lower and upper Boyd indices as follows,
\[
 \lbi(\sigma)=\lim_j \frac{\log_2 \lcst_j}{j}
 \quad\text{and}\quad
 \ubi(\sigma)=\lim_j \frac{\log_2 \ucst_j}{j}.
\]
Since $(\log \lcst_j)_{j\in\N}$ is a subadditive sequence, such limits always exist \cite{49.0047.01}. In this paper, $\sigma$ will always stand for an admissible sequence and $M$ for a natural number, possibly zero.

Starting from the definitions of the pointwise H\"older spaces $\Lambda^\alpha(x_0)$ (with $\alpha >0$) and the generalized uniform H\"older spaces $\Lambda^\sigma(\R^d)$ introduced in \cite{unifhold}, we are naturally led to the following definition.
\BDef\label{def:sp}
Let $x_0 \in \R^d$; a continuous function $f\in L^\infty_\loc(\R^d)$ belongs to $\Lambda^{\sigma,M}(x_0)$ if there exist $C,J>0$ such that
\[
 \inf_{P\in \P[M]} \|f-P\|_{L^\infty(x_0+B_j)} \le C \sigma_j,
\]
for any $j\ge J$.
\EDef
We trivially have the following alternative definition for $\Lambda^{\sigma,M}(x_0)$.
\BDef\label{rem:polj}
A function $f\in L^\infty_\loc(\R^d)$ belongs to $\Lambda^{\sigma,M}(x_0)$ if and only if there exist $C,J>0$ such that, for any $j\ge J$, there exists a polynomial $P_{j}\in \P[M]$ for which
\begin{equation}\label{eq:def2}
 \sup_{h\in B_j} |f(x_0+h)-P_{j}(x_0+h)| \le C \sigma_j.
\end{equation}
\EDef

Sometimes, we will also need to impose a slightly stronger condition than continuity to a function.
\BDef
A function $f$ is uniformly H\"older if and only if there exists $\epsilon>0$ such that $f\in \Lambda^\epsilon(\R^d)$ (here a function belonging to $\Lambda^\epsilon(\R^d)$ is necessarily continuous).
\EDef

\subsection{Independence of the polynomial from the scale}
It is important to remark that the polynomial occurring in inequality~(\ref{eq:def2}) is a function of the scale $j$. However, for the classical H\"older spaces, such polynomial is independent of $j$. Here, we look under which conditions the independence still holds in the generalized case, i.e.\ under which conditions $P_{j}=P \in \P[M]$ for any $j\ge J$.

We will need the following Markov inequality (see e.g.\ \cite{dit:92}): Let $p\in(0,\infty]$, $k \in\{1,\ldots,d\}$ and $S\subset \R^d$ be a bounded convex set with non-empty interior; one has
\[
 \| D_k P \|_{L^p(S)}\le C n^2 \|P\|_{L^p(S)},
\]
for any $P\in \P[n-1]$, where $C$ is a function of $S$ (but is independent of $P$ and $n$). If $x_0\in\R^d$, we thus have
\begin{equation}\label{eq:Markov}
 \| D_k P\|_{L^\infty(x_0+ r B)} \le \frac{C n^2}{r} \| P\|_{L^\infty(x_0+ r B)},
\end{equation}
for any $r>0$ and any $P\in\P[n-1]$, where $C$ is a constant (and does not depend on $P$, $n$ or $r$).
\BLem\label{lem:peano}
If $f\in \Lambda^{\sigma,M}(x_0)$ with $M<\lbi(\sigma^{-1})$, the sequence of polynomials $(P_{j})_j$ occurring in~(\ref{eq:def2}) satisfies
\[
 \|D^\beta P_{k} - D^\beta P_{j}\|_{L^\infty(x_0+B_k)} \le C 2^{j|\beta|} \sigma_j,
\]
for any multi-index $\beta$ such that $|\beta|\le M$ and $k\ge j\ge J$.

In particular, $(D^\beta P_j(x_0))_j$ is a Cauchy sequence for any multi-index $\beta$ such that $|\beta|\le M$.
\ELem
\BProof
Using the Markov inequality, we get
\begin{eqnarray*}
 \|D^\beta P_{j} - D^\beta P_{j+1}\|_{L^\infty(x_0+ B_{j+1})}
 &\le&
 C 2^{|\beta|j} \| P_{j} - P_{j+1}\|_{L^\infty(x_0+ B_{j+1})} \\
 &\le& C 2^{|\beta|j} (\| P_{j} - f\|_{L^\infty(x_0+ B_{j+1})} \\
 && + \| f-P_{j+1}\|_{L^\infty(x_0+ B_{j+1})}) \\
 &\le& C 2^{|\beta|j} (\sigma_j +\sigma_{j+1}) \\
 &\le& C 2^{|\beta|j} \sigma_j
\end{eqnarray*}
for any $\beta$ such that $|\beta|\le M$. Therefore, if $k$ satisfies $k\ge j\ge J$, one gets
\begin{eqnarray*}
 \| D^\beta (P_{j}- P_{k})\|_{L^\infty(x_0+ B_k)}
 &\le& \sum_{l=j}^{k-1} \| D^\beta (P_{l}- P_{l+1})\|_{L^\infty(x_0+ B_k)} \\
 &\le& \sum_{l=j}^{k-1} \| D^\beta (P_{l}- P_{l+1})\|_{L^\infty(x_0+ B_{l+1})} \\
 &\le& C \sum_{l=j}^{k-1} 2^{|\beta|l} \sigma_l \\
 &\le& C 2^{|\beta|j} \sigma_j,
\end{eqnarray*}
which is the desired result.
\EProof
\BLem
If $f\in \Lambda^{\sigma,M}(x_0)$ with $M<\lbi(\sigma^{-1})$ and $(P_{j})_j$ is a sequence of polynomials satisfying inequality~(\ref{eq:def2}), for any multi-index $\beta$ such that $|\beta|\le M$, the limit
\begin{equation}\label{eq:pean der}
 f_\beta(x_0)= \lim_{j\to\infty} D^\beta P_j(x_0)
\end{equation}
is independent of the chosen sequence $(P_{j})_j$.
\ELem
\BProof
If $(Q_{j})_j$ is another sequence of polynomials satisfying inequality~(\ref{eq:def2}), one gets
\[
 |D^\beta Q_{j}(x_0)- f_\beta(x_0)|\le
 |D^\beta Q_{j}(x_0)- D^\beta P_{j}(x_0)| + |D^\beta P_{j}(x_0) -f_\beta(x_0)|.
\]
Since one has, using the Markov inequality,
\begin{eqnarray*}
 \| D^\beta (P_{j} -Q_{j}) \|_{L^\infty(x_0+ B_j)}
 &\le&
 C 2^{|\beta|j} \| P_{j} -Q_{j} \|_{L^\infty(x_0+ B_j)} \\
 &\le&
 C 2^{|\beta|j} (\| P_{j} -f \|_{L^\infty(x_0+ B_j)} \\
 && +\| f -Q_{j} \|_{L^\infty(x_0+ B_j)}) \\
 &\le& C 2^{|\beta|j} \sigma_j\to 0,
\end{eqnarray*}
as $j\to\infty$, one can conclude.
\EProof
For such functions, we can introduce the notion of Peano derivative; this definition is similar to the ones given in \cite{dev:84,marianne}.
\BDef
Under the hypothesis of lemma~\ref{lem:peano}, the $\beta$-th Peano derivative of $f$ at $x_0$ is $f_\beta (x_0)= \lim_j D^\beta P_j(x_0)$.
\EDef

We can now obtain the result concerning the independence of the polynomials.
\BThe\label{thm:pol ind}
If $M<\lbi(\sigma^{-1})$, then $f\in \Lambda^{\sigma,M}(x_0)$ if and only if there exist $C>0$ and a unique polynomial $P\in \P[M]$ such that
\begin{equation}\label{eq:polind}
 \| f- P \|_{L^\infty(x_0+B_j)} \le C \sigma_j,
\end{equation}
for $j$ sufficiently large.
\EThe
\BProof
Let $(P_{j})_j$ be a sequence of polynomials for which inequality~(\ref{eq:def2}) is satisfied and set
\[
 P(x)= \sum_{|\beta|\le M} f_\beta(x_0) \frac{(x-x_0)^\beta}{|\beta|!}.
\]
One has
\begin{eqnarray*}
 \lefteqn{\|P- P_{j} \|_{L^\infty(x_0+B_j)}}&& \\
 &=& \| \sum_{|\beta|\le M} (f_\beta(x_0)- D^\beta P_{j}(x_0)) \frac{(x-x_0)^\beta}{|\beta|!} \|_{L^\infty(x_0+B_j)} \\
 &\le&
  \sum_{|\beta|\le M} |f_\beta(x_0)- D^\beta P_{j}(x_0)| 2^{-j|\beta|}.
\end{eqnarray*}

Since lemma~\ref{lem:peano} implies
\[
 |f_\beta(x_0) - D^\beta P_{j}(x_0)| \le C 2^{j|\beta|} \sigma_j,
\]
for $j$ sufficiently large, we have
\[
 \|P - P_{j} \|_{L^\infty(x_0+B_j)} \le C \sigma_j.
\]

This inequality can be used to obtain
\begin{eqnarray*}
  \|f- P \|_{L^\infty(x_0+B_j)}
 &\le&
  \|f- P_{j} \|_{L^\infty(x_0+B_j)}
 + \|P_{j}- P \|_{L^\infty(x_0+B_j)} \\
 &\le& C \sigma_j,
\end{eqnarray*}
which shows the existence of $P$.

If two polynomials $P,Q\in \P[M]$ satisfy inequality~(\ref{eq:polind}),
\[
 \|P -Q \|_{L^\infty(x_0+B_j)}
 \le \|P -f\|_{L^\infty(x_0+B_j)} +\| f-Q \|_{L^\infty(x_0+B_j)} \le C \sigma_j,
\]
but if $P\not=Q$,
\[
 \| P-Q \|_{L^\infty(x_0+B_j)}\ge C 2^{-j M},
\]
for $j$ sufficiently large, so that $2^{jM} \sigma_j$ does not tend to zero.
\EProof
The polynomial $P$ in inequality~(\ref{eq:polind}) is the Taylor expansion of $f$, where the derivative is replaced by the Peano derivative.

Let $\alpha\in (0,\infty)$; the sequence $\sigma=(2^{-j\alpha})_{j\in\N}$ is an admissible sequence with $\lbi(\sigma)=\ubi(\sigma)=-\alpha$, $\lbi(\sigma^{-1})=\ubi(\sigma^{-1})=\alpha$ and $\Lambda^\alpha(x_0)=\Lambda^{\sigma,\floor{\alpha}}(x_0)=\Lambda^{\sigma,[\lbi(\sigma^{-1})]}(x_0)$.
The definition given by (\ref{eq:us hol}) is very often slightly changed (we will use such a modified version in the sequel). It is easy to check that both definitions lead to the same spaces.
\BRem
It is easy to check that the polynomial satisfying equation (\ref{eq:us hol}) is unique if and only if $\alpha\not\in \N$. If $\alpha\in\N$, one rather imposes $P\in \P[\alpha-1]$ in order to obtain the uniqueness of the polynomial, so that $\Lambda^\alpha(x_0)= \Lambda^{\sigma, \lbi(\sigma^{-1})-1} (x_0)$, with $\sigma_j=2^{-j\alpha}$.
%
\ERem

The following proposition rigorously expresses the idea that the space $\Lambda^{\sigma,M}(x_0)$ associated to a sequence $(\sigma_j)_j$ that decreases faster than $2^{-jM}$ is included in the usual H\"older space $\Lambda^M(x_0)$.
\BCor
If $\lbi(\sigma^{-1})>M$, one has $\Lambda^{\sigma,M}(x_0) \subset \Lambda^M(x_0)$. 
\ECor
\BProof
Let $f\in\Lambda^{\sigma,M}(x_0)$, $P$ be defined as in theorem~\ref{thm:pol ind}, i.e.
\[
 P(x)= \sum_{|\beta|\le M} f_\beta(x_0) \frac{(x-x_0)^\beta}{|\beta|!}
\]
and let us set
\[
 Q(x)= \sum_{|\beta|\le M-1} f_\beta(x_0) \frac{(x-x_0)^\beta}{|\beta|!}.
\]

One gets
\begin{eqnarray*}
 \|f -Q \|_{L^\infty(x_0+B_j)}
 &\le& \| f-P \|_{L^\infty(x_0+B_j)}
 +  \| P-Q \|_{L^\infty(x_0+B_j)} \\
 &\le& C \sigma_j + C 2^{-jM} \le C 2^{-jM},
\end{eqnarray*}
since $2^{jM}\sigma_j$ tends to zero.
\EProof

\section{Alternative definitions of generalized H\"older spaces}
Since the uniform spaces $\Lambda^\sigma(\R^d)$ can be defined via finite differences or convolutions, one can wonder if such characterizations also hold for the pointwise version of these spaces.

\subsection{Characterization in terms of finite differences}
As usual, $\Delta_h^n f$ will stand for the finite difference of order $n$: given a function $f$ defined on $\R^d$ and $x,h\in\R^d$,
\[
 \Delta_h^1 f(x)=f(x+h)-f(x)
 \quad\text{and}\quad
 \Delta_h^{n+1} f(x)= \Delta_h^1 \Delta_h^n f(x),
\]
for any $n\in\N$. We also set
\[
 B_h^M(x_0,j)=\{x: [x,x+(M+1)h] \subset x_0+B_j\}
\]
In order to obtain a more general result, we drop the continuity condition of definition~\ref{def:sp} in this section.

\BPro
Let $f\in L^\infty_\loc(\R^d)$; one has $f\in \Lambda^{\sigma,M}(x_0)$ if and only if there exist $C,J>0$ such that
\begin{equation}\label{eq:diffin}
 \sup_{h\in B_j} \| \Delta_h^{M+1} f\|_{L^\infty(B_h^M(x_0,j))} \le C \sigma_j,
\end{equation}
for any $j\ge J$.
\EPro
\BProof
The theorem of Whitney (see e.g.\ \cite{bru70}) directly implies that if $f$ satisfies inequality~(\ref{eq:diffin}), then $f\in \Lambda^{\sigma,M}(x_0)$: One has
\[
 \inf_{P\in \P[M]} \| f-P\|_{L^\infty(x_0+B_j)} \le
 C \sup_{h\in B_j} \| \Delta_h^{M+1} f\|_{L^\infty(B_h^M(x_0,j))}.
\]

Let us now suppose that $f\in \Lambda^{\sigma,M}(x_0)$ and let $x\in B_h^M(x_0,j)$, $P\in \P[M]$. One has
\begin{eqnarray*}
 \|\Delta_h^{M+1} f\|_{L^\infty(B_h^M (x_0,j))} &=& \|\Delta_h^{M+1} (f-P)\|_{L^\infty(B_h^M (x_0,j))} \\
 &\le& (M+1)! (M+2) \|f-P\|_{L^\infty(x_0+ B_j)}.
\end{eqnarray*}

Now, since there exists a polynomial $P_{j}\in \P[M]$ such that
\[
 \|f -P_{j} \|_{L^\infty(x_0+B_j)} \le C \sigma_j,
\]
for $j$ sufficiently large, one gets
\[
 \sup_{h\in B_j} \| \Delta_h^{M+1} f\|_{L^\infty(B_h^M(x_0,j))}
 \le C \sigma_j,
\]
for $j$ sufficiently large.
\EProof

\subsection{Characterization in terms of convolutions}
Let us denote the space of the infinitely differentiable functions with compact support included in $E$ by $C^\infty_c(E)$. In this section, $\rho$ will denote a radial function of $C^\infty_c(B)$ such that $\rho(x)\in [0,1]$ for any $x\in\R^d$ and $\|\rho\|_{1}=1$. Moreover, one sets $\rho_j = 2^{-jd} \rho(\cdot/2^j)$, for any $j\in\N$.

In \cite{unifhold}, the following result has been obtained:
\BLem\label{lem:conv}
Let $N\in \N_0$; if $f\in L^1_\loc(\R^d)$ satisfies
\[
 \sup_{k\ge j} \| f* \rho_{k} -f \|_{L^\infty(x_0+B_j)} \le C \sigma_j,
\]
for $j\ge J$, then, for any multi-index $\beta$ such that $|\beta|\le N$, one has
\[
 \| D^\beta (f*\rho_{j} -f*\rho_{j-1}) \|_{L^\infty(x_0+B_j)} \le C 2^{jN} \sigma_j,
\]
for $j\ge J$.
\ELem

Using the same ideas as in \cite{unifhold}, one gets the desired characterization.
\BThe
If $f\in \Lambda^{\sigma,M}(x_0)$, then there exists a function $\Phi\in C^\infty_c(\R^d)$ such that
\begin{equation}\label{eq:holconv}
 \sup_{k\ge j} \| f- f*\Phi_k\|_{L^\infty(x_0 +B_j)} \le C\sigma_j,
\end{equation}
for $j$ sufficiently large.

Conversely, if $\sigma\to 0$, $f$ is uniformly H\"older and if $f$ satisfies inequality~(\ref{eq:holconv}) for a function $\Phi\in C^\infty_c(\R^d)$, then $f\in \Lambda^{\sigma,M}(x_0)$ for any $M\in\N_0$ such that $M+1>\ubi(\sigma^{-1})$.
\EThe
\BProof
As in \cite{unifhold} (see also \cite{kra:83}), let us set
\[
 \Psi(x)= \sum_{j=0}^{m/2-1} (-1)^j \frac{m!}{j! (m-j)!} \frac{1}{2j-m} \rho(\frac{x}{2j-m}),
\]
where $m$ is large enough (larger than $M+1$) and $\Phi= \Psi/\int \Psi dx$. Using the same arguments as in \cite{unifhold}, one gets
\[
 f*\Phi_k(x) -f(x)= C \int \Delta_{2^{-k}t}^m f(x) \rho(t) \,dt,
\]
which leads to inequality~(\ref{eq:holconv}).

Let us show the converse. Let $\alpha\in (0,1)$ such that $f\in\Lambda^\alpha(\R^d)$ and set, as in \cite{unifhold},
\[
 f_1=f* \Phi_1
 \quad\text{and}\quad
 f_j= f*(\Phi_j -\Phi_{j-1}),
\]
for $j>1$. Since $f$ is uniformly H\"older, $f$ is uniformly equal to $\sum_{j\ge 1} f_j$ on $\R^d$ and
\[
 \Delta_{h}^{M+1} f= \sum_{j\ge 1} \Delta_h^{M+1} f_j
\]
uniformly on $\R^d$, for any $h\in \R^d$. For $j\in\N$, let $n_0\in\N_0$, $h\in\R^d$ and $j_0\in\N_0$ be such that $M+1<2^{n_0}$, $ |h|\le 2^{-(j+n_0)}$ and $2^{-(j_0+1)\alpha}\le \sigma_j \le 2^{-j_0\alpha}$. One has
\begin{eqnarray*}
 \|\Delta_h^{M+1} f\|_{L^\infty(x_0+B_j)} &\le&
 \sum_{k=1}^{j-1} \|\Delta_h^{M+1} f_k\|_{L^\infty(x_0+B_j)}
 +\sum_{k=j}^{j_0} \|\Delta_h^{M+1} f_k\|_{L^\infty(x_0+B_j)} \\
 && + \sum_{k\ge j_0+1} \|\Delta_h^{M+1} f_k\|_{L^\infty(x_0+B_j)},
\end{eqnarray*}
where the second term in the majoration only appears if $j\le j_0$.

Using lemma~\ref{lem:conv} and the fact that $M+1>\ubi(\sigma^{-1})$, the mean value theorem allows to write
\begin{eqnarray*}
 \sum_{k=1}^{j-1} \|\Delta_h^{M+1} f\|_{L^\infty(x_0+B_j)} &\le&
 \sum_{k=1}^{j-1} C |h|^{M+1} \sup_{|\beta|=M+1} \| D^\beta f_k\|_{L^\infty(x_0+B_{j-1})} \\
 &\le& C 2^{-j(M+1)} \sum_{k=1}^{j-1} 2^{k(M+1)} \sigma_k
 \le C \sigma_j.
\end{eqnarray*}
Moreover,
\begin{eqnarray*}
 \sum_{k=j}^{j_0} \|\Delta_h^{M+1} f_k\|_{L^\infty(x_0+B_j)} &=&
 \| \Delta_h^{M+1} (f*\Phi_{j_0} -f* \Phi_{j-1}) \|_{L^\infty(x_0+B_j)} \\
 &\le& C \| f*\Phi_{j_0} -f* \Phi_{j-1} \|_{L^\infty(x_0+B_{j-1})} \\
 &\le& C (\| f*\Phi_{j_0} -f\|_{L^\infty(x_0+B_{j-1})} \\
 && + \|f -f* \Phi_{j-1} \|_{L^\infty(x_0+B_{j-1})}) \\
 &\le& C 2^{-j_0\alpha} + C \sigma_{j-1} \le C \sigma_j.
\end{eqnarray*}
Finally,
\begin{eqnarray*}
 \sum_{k\ge j_0+1} \|\Delta_h^{M+1} f_k\|_{L^\infty(x_0+B_j)}
 &\le& C \sum_{k\ge j_0+1} \|f_k\|_{L^\infty(\R^d)} \\
 &\le&C \sum_{k\ge j_0+1} 2^{-k\alpha} \\
 &\le& 2^{-j_0\alpha} \le C \sigma_j.
\end{eqnarray*}
One then has,
\begin{eqnarray*}
 \sup_{h\in B_{j+n_0}} \| \Delta_h^{M+1} f\|_{L^\infty(x_0+B_{j+n_0})}
 &\le& \sup_{h\in B_{j+n_0}} \| \Delta_h^{M+1} f\|_{L^\infty(x_0+B_j)} \\
 &\le& C\sigma_j\le C \sigma_{j+n_0},
\end{eqnarray*}
as wanted.
\EProof

\section{Generalized pointwise H\"older spaces and wavelets}
The usual H\"older spaces can ``nearly'' be characterized in terms of wavelet \cite{jaffard}: for the sufficiency of the condition, the function has to be uniformly H\"older and a logarithmic correction appears. We show here that such a result still holds in the generalized case.

\subsection{Definitions}
Let us briefly recall some definitions and notations (for more precisions, see e.g.~\cite{dau:92,mey:90,mal:98}). Under some general assumptions, there exist a function $\phi$ and $2^d-1$ functions $(\psi^{(i)})_{1\le i<2^d}$, called wavelets, such that
\[
 \{\phi(\cdot -k):k\in\Z^d \}
 \cup
 \{\psi^{(i)}(2^j \cdot -k):1\le i<2^d, k\in \Z^d, j\in\N_0 \}
\]
form an orthogonal basis of
$L^2(\R^d)$. Any function $f\in L^2(\R^d)$ can be decomposed as
follows,
\[
f(x)=\sum_{k\in \Z^d} C_k \phi(x-k) + \sum_{j=0}^{+\infty}
\sum_{k\in\Z^d} \sum_{1\le i<2^d} c^{(i)}_{j,k} \psi^{(i)}(2^j x-k),
\]
where
\[
c^{(i)}_{j,k}=2^{dj}\int_{\R^d}f(x) \psi^{(i)}(2^jx-k)\, dx,
\]
and
\[
C_k=\int_{\R^d} f(x) \phi(x-k)\, dx.
\]
Let us remark that we do not choose the $L^2(\R^d)$ normalization for
the wavelets, but rather a $L^\infty$ normalization, which is
better fitted to the study of the H\"{o}lderian regularity. Hereafter, the
wavelets are always supposed to belong to $C^n(\R^d)$ with
$n>M$, and the functions
$(D^{s}\phi)_{|s|\le \gamma}$,
$(D^{s}\psi^{(i)})_{|s|\le \gamma}$ are assumed to have
fast decay.

A dyadic cube of scale $j$ is a cube of the form
\[
\lambda=\left[\frac{k_1}{2^j},\frac{k_1+1}{2^j}\right)\times \cdots
\times \left[\frac{k_d}{2^j},\frac{k_d+1}{2^j}\right),
\]
where $k=(k_1,\ldots,k_d)\in \Z^d$.
From now on, wavelets and wavelet coefficients will be indexed
with dyadic cubes $\lambda$. Since $i$ takes $2^d-1$ values, we can assume
that it takes values in $\{0,1\}^d-(0,\ldots,0)$; we will use the
following notations:
\begin{itemize}
\item $\displaystyle \lambda=\lambda(i,j,k)=\frac{k}{2^j}+\frac{i}{2^{j+1}}+\left[0,\frac{1}{2^{j+1}}\right)^d$,
\item $c_\lambda=c^{(i)}_{j,k}$,
\item $\psi_\lambda=\psi^{(i)}_{j,k}=\psi^{(i)}(2^j\cdot -k)$.
\end{itemize}
The pointwise H\"olderian regularity of a function is closely related
to the decay rate of its wavelet leaders.
\BDef
The wavelet leaders are defined by
\[
d_\lambda=\sup_{\lambda'\subset\lambda} |c_{\lambda'}|.
\]
Two dyadic cubes $\lambda$ and $\lambda'$ are adjacent if they are
at the same scale and if $\dist(\lambda,\lambda')=0$. We denote by
$3\lambda$ the set of $3^d$ dyadic cubes adjacent to $\lambda$ and
by $\lambda_j(x_0)$ the dyadic cube of side of length $2^{-j}$ containing
$x_0$; then
\[
d_j(x_0)=\sup_{\lambda\subset 3\lambda_j(x_0)} d_{\lambda}.
\]
\EDef

\subsection{Result}
From now on, we will suppose that the wavelets are compactly supported; such wavelets are constructed in \cite{dau:88} and $j_0$ will stand for a natural number such that the support of $\psi^{(i)}$ is included in $2^{j_0}B$, for any $i\in \{1,\ldots, 2^d-1\}$.

\BThe
If $f\in \Lambda^{\sigma,M}(x_0)$, there exist $C>0$ and $J\in\N_0$ such that
\begin{equation}\label{eq:thm ondel}
 d_j(x_0) \le C\sigma_j,
\end{equation}
for any $j\ge J$.

Conversely, let $f$ be an uniformly H\"older function; if inequality~(\ref{eq:thm ondel}) is satisfied for an admissible sequence $\sigma$ that tends to zero, then $f\in \Lambda^{\tau,M} (x_0)$, where $\tau$ is the admissible sequence defined by $\tau_j=\sigma_j |\log_2 \sigma_j|$ and $M\in\N_0$ is any number satisfying $M+1>\ubi(\sigma^{-1})$.
\EThe
\BProof
If $f\in \Lambda^{\sigma,M} (x_0)$, let $k_0\in \N_0$ be such that $2^{j_0+1}+4d \le 2^{k_0}$. For $j\ge k_0+1$ and $\lambda=\lambda(i,j',k')\subset 3\lambda_j(x_0)$, one has
\begin{eqnarray*}
 |c_\lambda| &=&
 | 2^{dj'} \int f(x) \psi_\lambda(x) \,dx| \\
 &=& | 2^{dj'} \int \big(f(x)-P_{j-k_0}(x)\big) \psi_\lambda(x) \, dx| \\
 &=& | 2^{dj'} \int_{\frac{k'}{2^{j'}}+ B_{j'-j_0}} \big(f(x)-P_{j-k_0}(x)\big) \psi_\lambda(x) \, dx| \\
 &\le& 2^{dj'} \int_{x_0+ B_{j-k_0}} \big| f(x)-P_{j-k_0}(x)\big|\ |\psi_\lambda(x)| \, dx \\
 &\le& C 2^{dj'} \sigma_{j-k_0} \int |\psi_\lambda(x)| \,dx \le C \sigma_j,
\end{eqnarray*}
which is the desired result.

Now, let us suppose that inequality~(\ref{eq:thm ondel}) is satisfied for a function $f\in \Lambda^\epsilon(\R^d)$. Let us set
\[
 f_{-1}= \sum_k C_k \phi(\cdot -k)
 \quad\text{and}\quad
 f_j=\sum_{i,k} c_\lambda \psi_\lambda,
\]
for $j\in\N_0$. In \cite{unifhold2}, it has been shown that these functions have the same regularity as the wavelets and that $f$ is uniformly equal to $\sum_{j\ge -1} f_j$. Let us define
\[
 P_{J}(x)= \sum_{|\beta|\le M} \frac{(x-x_0)^\beta}{|\beta|!} \sum_{j=-1}^J D^\beta f_j(x_0)
\]
and let us choose $n_d\in\N$ such that  $R>2^{-j}$ and $k/2^j\in x+RB$ ($x\in\R^d$) implies
\[
 \frac{k}{2^j}+\frac{i}{2^{j+1}}+[0,\frac{1}{2^j})^d\subset x +2^{n_d} RB.
\]
Let us also choose $m_d\in\N$ such that any ball $x+ B_j$ ($x\in\R^d$, $j\in\N_0$) is included in a dyadic cube of length $2^{m_d-j}$. If $J'$ is such that $\sigma_j<1$ for any $j\ge J'$, we finally choose $J$ such that $J\ge \{J', j_0+n_d+m_d+1\}$. One has
\begin{eqnarray*}
 \lefteqn{\|f-P_{J} \|_{L^\infty(x_0+ B_J)}} \\
 &\le& \sum_{j=-1}^J \| f_j(x) -\sum_{|\beta|\le M} \frac{(x-x_0)^\beta}{|\beta|!} D^\beta f_j(x_0) \|_{L^\infty(x_0+ B_J)} \\
 && + \sum_{j\ge J+1} \| f_j\|_{L^\infty(x_0+ B_J)}
\end{eqnarray*}

Let us look at the first term of the majoration. Let $j\le J$; using the Taylor expansion, one gets
\begin{eqnarray*}
 \lefteqn{\| f_j(x) -\sum_{|\beta|\le M} \frac{(x-x_0)^\beta}{|\beta|!} D^\beta f_j(x_0) \|_{L^\infty(x_0+ B_J)}} \\
 &\le& C 2^{-J(M+1)} \sup_{|\beta|=M+1} \| D^\beta f_j\| _{L^\infty(x_0+ B_J)}.
\end{eqnarray*}
If $\beta$ satisfies $|\beta|=M+1$, we have, for any $x\in x_0+ B_J$,
\begin{eqnarray*}
 |D^\beta f_j (x)| &\le& \sum_{i,k} 2^{j(M+1)} |c_\lambda| |D^\beta \psi_\lambda(x)| \\
 &=& \sum_i \sum_{k2^{-j}\in x+ B_{j-j_0}} 2^{j(M+1)} |c_\lambda| |D^\beta \psi_\lambda(x)|.
\end{eqnarray*}
Each coefficient $c_\lambda$ in the last sum is such that $\lambda\subset x+ B_{j-j_0-n_d}$. Therefore, if $j\ge j_0+n_d+m_d+1$,
\[
 |c_\lambda|\le C \sigma_{j+j_0+n_d+m_d+1}.
\]
Otherwise, since $f$ is uniformly H\"older, $|c_\lambda|\le C \le C\sigma_j$. Therefore,
\[
 \|D^\beta f_j\|_{L^\infty(x_0+ B_J)}
 \le C 2^{j(M+1)} \sigma_j,
\]
which implies
\begin{eqnarray*}
  \lefteqn{\| f_j(x) -\sum_{|\beta|\le M} \frac{(x-x_0)^\beta}{|\beta|!} D^\beta f_j(x_0) \|_{L^\infty(x_0+ B_J)}} \\
 &\le& C 2^{-J(M+1)} \sum_{j=-1}^J 2^{j(M+1)} \sigma_j \le C\sigma_J.
\end{eqnarray*}

For the second term in the majoration, let us define $J_1\in\N$ as the number such that $2^{-\epsilon J_1}\le \sigma_J < 2^{-\epsilon(J_1-1)}$ and decompose the sum as follows:
\[
 \sum_{j\ge J+1} \| f_j\|_{L^\infty(x_0+ B_J)} =  \sum_{j\ge J_1+1} \| f_j\|_{L^\infty(x_0+ B_J)} + \sum_{j= J+1}^{J_1} \| f_j\|_{L^\infty(x_0+ B_J)}
\]
We have
\begin{eqnarray*}
 \sum_{j\ge J_1+1} \| f_j\|_{L^\infty(x_0+ B_J)} &\le&
 \sum_{j\ge J_1+1} \|f_j\|_{L^\infty(\R^d)}
 \le C \sum_{j\ge J_1+1} 2^{-\epsilon j} \\
 &\le& C 2^{-\epsilon J_1}
 \le C \sigma_J.
\end{eqnarray*}
Now, for $j\in \{J+1,\ldots,J_1\}$ and $x\in x_0+ B_J$, one has
\[
 |f_j(x)| \le \sum_i \sum_{k2^{-j} \in x+ B_{j-j_0}} |c_\lambda \psi_\lambda(x)|.
\]
If $j\ge J+j_0+n_d$, the wavelet coefficients $c_\lambda$ in the last sum are such that
\[
 \lambda\subset x+ B_{j-j_0-n_d} \subset x_0+B_{J-1}
\]
and therefore
\[
 |c_\lambda|\le C \sigma_{J-m_d-1}\le C \sigma_J.
\]
In the other case,
\[
 \lambda\subset x+ B_{j-j_0-n_d}\subset x_0 + B_{j-j_0-n_d-1}
\]
and thus
\[
 |c_\lambda|\le C \sigma_{j-j_0-n_d-m_d-1}\le C\sigma_j\le C\sigma_J.
\]
These inequalities lead to
\[
 \sum_{j= J+1}^{J_1} \| f_j\|_{L^\infty(x_0+ B_J)} \le C J_1 \sigma_J \le C |\log_2(\sigma_J)| \sigma_J.
\]
Putting all these inequalities together, one gets
\[
 \|f- P_{J}\|_{L^\infty(x_0+ B_J)} \le C |\log_2(\sigma_J)| \sigma_J,
\]
as desired.
\EProof

The converse part of the previous theorem requires a uniform regularity condition. As shown in \cite{jaf:91}, a stronger condition than continuity is necessary in the usual case (see also \cite{jaf:00}, where similar results are obtained (in the usual case) with a Besov regularity assumption). Similarly, the logarithmic correction is best possible in the usual case \cite{jaf:91}.

\section{A generalized definition of the H\"older exponent}
The usual H\"older spaces are embedded: $\alpha<\beta$ implies $\Lambda^\beta(x_0)\subset \Lambda^\alpha(x_0)$. A notion of regularity for a function $f\in L^\infty_\loc(\R^d)$ at $x_0$ can thus be given by the so-called H\"older exponent,
\[
 h_f(x_0)= \sup \{ \alpha>0: f\in \Lambda^\alpha (x_0) \}.
\]
To do so in the generalized case, one needs some conditions under which $\Lambda^{\sigma,M}(x_0) \subset \Lambda^{\sigma',M'}(x_0)$.

\subsection{Preliminary results}
We first need some technical easy results. 
%
From now on, if $f\in \Lambda^{\sigma,M}(x_0)$, $(P_{j})_j$ will stand for the sequence of polynomials of $\P[M]$ corresponding to the definition. 
We will write
\[
 P_{j}(x)=\sum_{|\beta|\le M} a_j^{(\beta)} x^\beta
\]
and
\[
 Q_{j}(x)=\sum_{|\beta|\le M-1} a_j^{(\beta)} x^\beta
\]

\BLem\label{lem:fam hol}
Let $f\in \Lambda^{\sigma,M}(x_0)$; one has
\begin{equation}\label{eq:lem hol1}
 \sup_{|\beta|=M} |a_j^{(\beta)}| \le C (\sum_{k=1}^{j-1} (2^M \ucst_1)^k +1)
\end{equation}
and
\begin{equation}\label{eq:lem hol2}
 \sup_{|\beta|=M} |a_j^{(\beta)}| \le C (\sigma_j \lcst_1^{-j} \sum_{k=1}^{j-1} (2^M \lcst_1)^k +1).
\end{equation}
\ELem
\BProof
Using the Markov inequality~(\ref{eq:Markov}), we get
\begin{eqnarray*}
 \lefteqn{\| D^\beta (P_{j} - P_{j+1})\|_{L^\infty(x_0+ B_{j+1})}} \\
 &\le& C 2^{jM} \| P_{j} - P_{j+1}\|_{L^\infty(x_0+ B_{j+1})} \\
 &\le& C 2^{jM} (\| P_{j} - f\|_{L^\infty(x_0+ B_j)} +\| f - P_{j+1}\|_{L^\infty(x_0+ B_{j+1})}) \\
 &\le& C 2^{jM} \sigma_j,
\end{eqnarray*}
for any $\beta$ such that $|\beta|\le M$ and $j$ sufficiently large.
Therefore, we have
\begin{eqnarray*}
 \lefteqn{\| D^\beta (P_{1} - P_{j})\|_{L^\infty(x_0+ B_j)}} \\
 &\le& \sum_{k=1}^{j-1} \| D^\beta (P_{k} - P_{k+1})\|_{L^\infty(x_0+ B_j)} \\
 &\le& \sum_{k=1}^{j-1} \| D^\beta (P_{k} - P_{k+1})\|_{L^\infty(x_0+ B_{k+1})} \\
 &\le& C \sum_{k=1}^{j-1} 2^{kM} \sigma_k\le C \sum_{k=1}^{j-1} (2^{M} \ucst_1)^k,
 \end{eqnarray*}
for any $j$.

Now, let $\beta$ be a multi-index such that $|\beta|=M$; inequality~(\ref{eq:lem hol1}) follows from
\[
 \| D^\beta (P_{1} - P_{j})\|_{L^\infty(x_0+ B_j)} \ge C(|a_j^{(\beta)}|- |a_1^{(\beta)}|),
\]
while inequality~(\ref{eq:lem hol2}) can be obtained in the same way, using
\[
 \| D^\beta (P_{1} - P_{j})\|_{L^\infty(x_0+ B_j)} \le C \sum_{k=1}^{j-1} 2^{kM} \sigma_k
 \le C \sigma_j \lcst_1^{-j} \sum_{k=1}^{j-1} (2^M \lcst_1)^k,
\]
valid for any $j$.
\EProof

\BCor
Let $f\in\Lambda^{\sigma,M}(x_0)$; we have the following inequalities:
\begin{itemize}
 \item if $2^M \ucst_1<1$,
\[
 \|f- Q_{j}\|_{L^\infty(x_0+ B_j)} \le C (\sigma_j +2^{-jM}),
\]
 \item if $2^M \ucst_1>1$,
\[
 \|f- Q_{j}\|_{L^\infty(x_0+ B_j)} \le C (\sigma_j + \ucst_1^j),
\]
 \item if $2^M \ucst_1=1$,
\[
 \|f- Q_{j}\|_{L^\infty(x_0+ B_j)} \le C (\sigma_j + 2^{-jM}j).
\]
\end{itemize}
\ECor
\BCor
Let $f\in\Lambda^{\sigma,M}(x_0)$; we have the following inequalities:
\begin{itemize}
 \item if $2^M \lcst_1<1$,
\[
 \|f- Q_{j}\|_{L^\infty(x_0+ B_j)} \le C (\sigma_j (2^M \lcst_1)^{-j}+2^{-jM}),
\]
 \item if $2^M \lcst_1>1$,
\[
 \|f- Q_{j}\|_{L^\infty(x_0+ B_j)} \le C (\sigma_j + 2^{-jM}),
\]
 \item if $2^M \lcst_1=1$,
\[
 \|f- Q_{j}\|_{L^\infty(x_0+ B_j)} \le C (\sigma_j j + 2^{-jM}).
\]
\end{itemize}
\ECor

\subsection{Definitions}
Before introducing a definition of generalized H\"older exponent, we must first consider embedded spaces of type $\Lambda^{\sigma,M}(x_0)$. Once the definitions given, we provide sufficient conditions for generalized H\"older spaces to be embedded.

\BDef
If for any $\alpha>0$, $\sigma^{(\alpha)}$ is an admissible sequence, the application
\[
 \sigma^{(\cdot)}: \alpha>0 \mapsto \sigma^{(\alpha)}
\]
is called a family of admissible sequences.
\EDef
\BDef
Let $x_0\in\R^d$; a family $\sigma^{(\cdot)}$ of admissible sequences is decreasing for $x_0$ if $\alpha<\beta$ implies $\Lambda^{\sigma^{(\alpha)},\floor{\alpha}}(x_0) \subset \Lambda^{\sigma^{(\beta)},\floor{\beta}}(x_0)$.
\EDef
\BDef
 Let $\sigma^{(\cdot)}$ be a decreasing family of admissible sequences for $x_0$; if $f \in L^\infty_\loc(\R^d)$, the H\"older exponent of $f$ at $x_0$ for the family $\sigma^{(\cdot)}$ is given by
\[
 h^{\sigma^{(\cdot)}}_f (x_0)= \sup \{ \alpha>0: f\in \Lambda^{\sigma^{(\alpha)},\floor{\alpha}} (x_0) \}.
\]
\EDef

The following proposition is a simple corollary of the results obtained in the previous section; it helps to check if a family of admissible sequences is decreasing. If $\sigma^{(\cdot)}$ is a family of admissible sequences, we set
\[
 \lcst_j^{(\alpha)}=\inf_{k\in\N} \frac{\sigma_{j+k}^{(\alpha)}}{\sigma_k^{(\alpha)}}
 \quad\text{and}\quad
 \ucst_j^{(\alpha)}=\sup_{k\in\N} \frac{\sigma_{j+k}^{(\alpha)}}{\sigma_k^{(\alpha)}}.
\]
\BPro\label{pro:dec fam}
Let $\sigma^{(\cdot)}$ be a family of admissible sequences and $x_0\in\R^d$; $\sigma^{(\cdot)}$ is decreasing for $x_0$ if it satisfies the two following conditions:
\begin{enumerate}
 \item if $m\le \alpha <\beta<m+1$, with $m\in\N_0$, there exist $C,J>0$ such that
\[
 \sigma^{(\beta)}_j \le C \sigma_j^{(\alpha)},
\]
for any $j\ge J$,

 \item for any $m\in\N$, at least one of the two following conditions is satisfied:
 \begin{enumerate}
  \item there exists $\epsilon_0>0$ such that, for any $\epsilon\in (0,\epsilon_0)$, there exist $C,J>0$ for which $\sigma^{(m)}_j\le C \sigma^{(m-\epsilon)}_j$ and
        \begin{itemize}
         \item if $1<2^m \ucst^{(m)}_1$, $(\ucst^{(m)}_1)^j \le C \sigma^{(m-\epsilon)}_j$,
         \item if $1>2^m \ucst^{(m)}_1$, $2^{-jm} \le C \sigma^{(m-\epsilon)}_j$,
         \item if $1=2^m \ucst^{(m)}_1$, $j2^{-jm} \le C \sigma^{(m-\epsilon)}_j$,
        \end{itemize}
        for any $j\ge J$,
  \item there exists $\epsilon_0>0$ such that, for any $\epsilon\in (0,\epsilon_0)$, there exist $C,J>0$ for which $2^{-jm} \le C \sigma^{(m-\epsilon)}_j$ and
        \begin{itemize}
         \item if $1<2^m \lcst^{(m)}_1$, $\sigma^{(m)}_j \le C \sigma^{(m-\epsilon)}_j$,
         \item if $1>2^m \lcst^{(m)}_1$, $\sigma^{(m)}_j (2^{m} \lcst^{(m)}_1)^{-j} \le C \sigma^{(m-\epsilon)}_j$,
         \item if $1=2^m \lcst^{(m)}_1$, $j \sigma^{(m)}_j \le C \sigma^{(m-\epsilon)}_j$,
        \end{itemize}
        for any $j\ge J$.
 \end{enumerate}
\end{enumerate}
\EPro
This result is similar to the one obtained in \cite{unifhold} (under the hypothesis of proposition~\ref{pro:dec fam}, one gets a decreasing family of admissible sequences for the uniform case), but the proof given for these generalized uniform H\"older spaces cannot be adapted for the pointwise case.


\end{document}